\documentclass[sts,preprint]{imsart}
\setattribute{journal}{name}{}
\usepackage{graphicx}
\usepackage{amsmath}
\usepackage{amsthm}
\usepackage{amsfonts}
\usepackage{amssymb}
\usepackage{color}
\usepackage{verbatim}

\theoremstyle{plain}
\newtheorem{thm}{Theorem}[section]

\newtheorem{prop}[thm]{Proposition}

\newtheorem*{thm*}{Theorem}

\theoremstyle{remark}

\newtheorem{example}[thm]{Example}

\theoremstyle{definition}
\newtheorem{defn}[thm]{Definition}
\newtheorem{ques}[thm]{Question}

\numberwithin{equation}{section}

\def\EE{\mathbb{E}}
\def\PP{\mathbb{P}}
\def\R{\mathbb{R}}
\def\N{\mathbb{N}}
\def\Z{\mathbb{Z}}

\def\cA{\mathcal{A}}

\def\cF{\mathcal{F}}

\def\cN{\mathcal{N}}

\def\cX{\mathcal{X}}
\def\cY{\mathcal{Y}}

\DeclareMathOperator{\argmax}{argmax}
\DeclareMathOperator{\argmin}{argmin}

\DeclareMathOperator{\Meas}{Meas}

\begin{document}

\begin{frontmatter}
\title{Statistical Inference for Dynamical Systems: A review}
\begin{aug}
  \author{Kevin McGoff\ead[label=e1]{mcgoff@math.duke.edu}},
  \author{Sayan Mukherjee\ead[label=e2]{sayan@stat.duke.edu}},
  \author{Natesh S. Pillai\ead[label=e3]{pillai@fas.harvard.edu}}
  
    \address{Department of Mathematics,\\
    Duke University,\\
    Durham, NC\\ USA, 27708\\
    \printead{e1}}
    \address{Departments of Statistical Science, Computer Science, and Mathematics,\\
    Duke University,\\
    Durham, NC\\ USA, 27708\\
    \printead{e2}}
   \address{Department of Statistics,\\
    Harvard University,\\
    1 Oxford Street, Cambridge\\  MA, USA, 02138\\
    \printead{e3}}
 
 \end{aug}
\maketitle

\begin{abstract}
The topic of statistical inference for dynamical systems has been studied extensively across several fields. In this survey we focus on the problem of parameter estimation for nonlinear dynamical systems. Our objective is to place results across distinct disciplines in a common setting and  highlight opportunities for further research.
\end{abstract}
\end{frontmatter}
\maketitle

%%%%%%%%%%%%%%%%%%%%%%%%%%%%%%%%%%%%%%%%%%%%%%%%%%%%%%%%%%%%%%%%%%%%
%%%%%%%%%%%%%%%%%%%%%%%%%%%%%%%%%%%%%%%%%%%%%%%%%%%%%%%%%%%%%%%%%%%%
\section{Introduction}

The problem of parameter estimation in dynamical systems appears in many areas of science and engineering. Often the form of the model can be derived from some knowledge about the process under investigation, but parameters of the model must be inferred from empirical observations in the form of time series data. As this problem has appeared in many different contexts, partial solutions to this problem have been proposed in a wide variety of disciplines, including nonlinear dynamics in physics, control theory in engineering, state space modeling in statistics and econometrics, and ergodic theory and dynamical systems in mathematics. One purpose of this study is to present these various approaches in a common language, with the hope of unifying some ideas and pointing towards interesting avenues for further study.

By a dynamical system we mean a stochastic process of the form $(X_n,Y_n)_n$, where $X_{n+1}$ depends only on $X_n$ and possibly some noise, and $Y_n$ depends only on $X_n$ and possibly some noise. We think of $X_n$ as the true state of the system at time $n$ and $Y_n$ as our observation of the system at time $n$. The case when no noise is present has been most often considered by mathematicians in the field of dynamical systems and ergodic theory. In this case, all uncertainty in the system comes from the uncertainty in the initial state of the system, and the ability to estimate any parameters in the system may depend strongly on properties of the observation function $f(X_n) = Y_n$, although such questions have rarely been addressed rigorously. State space models, considered most often by statisticians, lie at the other end of the noise spectrum, where both $X_{n+1}$ and $Y_n$ depend on some noise. Hidden Markov models, which have received considerable attention, provide a broad class of examples of these systems. In this setting, the statistical question of consistency for methods of parameter estimation has been studied, and some general results are available. The other two possible assumptions on the presence of noise (assuming only dynamical noise or assuming only observational noise) have received relatively little attention, especially from the statistical point of view. While many of the proposed methods of parameter estimation for dynamical systems with observational noise have been studied via numerical simulations or on particular data sets, very few of these methods have been studied on a theoretical level. In fact, in the observational noise setting, basic statistical questions, such as whether a proposed method is consistent, have been considered only rarely, if at all, despite the fact that such systems capture important features of many experimental settings.

Consider, for example, the question of parameter inference for models of gene regulatory networks. The underlying model often favored by biologists consists of a system of ordinary differential equations, with each variable in the state vector representing the expression level of a particular gene in the network. For some networks of interest, a significant amount of work has produced biological understanding regarding the qualitative interactions between the genes in the network, but the corresponding ODE models still contain several parameters necessary for quantifying these interactions. Experimentalists are able to conduct experiments in which the expression levels of the genes in the network are measured at regularly spaced instances of time. The resulting data may be interpreted as time series data generated by a system of ODEs with noisy observations. The parameter inference problem in this setting consists of inferring the parameters of the ODE model from the observed data, and to the best of our knowledge there are no general statistical inference schemes for this type of problem that have been shown to be consistent.

Another example of interest is identifying the behavior of a dynamical system on a network. In a variety of applications one considers nodes in a communication network and measures the states of these nodes (or properties of the nodes) over time. In many settings, one would like to detect drastic changes in the nature of the dynamic behavior of the system. This problem is of vital importance to a variety of security applications on networks, and it can be formalized as the inference of large changes in the parameters of the network -- a change point model for a dynamic network.

The objective of this article is to survey methodology across a variety of fields for parameter inference in stochastic dynamical systems. We first state the various goals of inference in dynamical systems. Our focus will be parameter inference and we provide a natural decomposition of parameter inference into four possible settings defined by the structure of noise in the system. We then state what is known in terms of rigorous results for parameter inference in these four settings. Of these settings the case of deterministic dynamics with observational noise is the least developed in terms of sound statistical theory and will be our focus. We also mention several important open problems for parameter inference in these types of systems.

There is an extremely large body of work stretching across many disciplines that relates to the topic of statistical properties of dynamical systems. Although we attempt to provide references when possible, we make no attempt to be exhaustive, and we recognize that in fact many references have been omitted. On the other hand, we hope that the references cited in this article may serve as an appropriate starting point for further reading.

\section{Basic definitions and Preliminaries}

The most general setting that we will consider may be described as follows. Let $\cA$, $\cX$, $\cY$, and $\cN$ be Polish \footnote{This is a classical assumption in dynamical systems.} spaces (complete metric spaces with a countable dense set), where each one is equipped with its Borel $\sigma$-algebra. The space $\cA$ denotes the parameter space, the underlying dynamical system evolves in the space $\cX$ and the observations take values in $\cY$. %The state space dynamics is given by a measurable map $Q : \cA \times \cX \times \cN \to \cX$, and the observations are governed by a measurable map $f : \cA \times \cX \times \cN \to \cY$. 
We consider a stochastic process $(X_n,Y_n)_n$, which satisfies the following dynamics: for some $a$ in $\cA$ and $X_0$ distributed according to a Borel probability measure on $\cX$,
\begin{align}
Y_{n} & = f_a(X_{n},\epsilon_{n}), \label{generalobs} \\
X_{n+1} & = T_a(X_n,\delta_{n+1}), \label{generaldyn}
\end{align}
where $\delta_{n+1}$ is the dynamical noise and $\epsilon_{n}$ is the observational noise. The maps
$T : \cA \times \cX \times \cN \to \cX$ and $f : \cA \times \cX \times \cN \to \cY$ determine the evolution of the state space dynamics and the observation process, respectively. We refer to a sequence $(X_n)_n$ satisfying (\ref{generaldyn}) as a trajectory and a sequence $(Y_n)_n$ satisfying (\ref{generalobs}) as a sequence of observations.

Let us start with the following definitions.
\begin{defn}
 A stochastic process $(X_n)_n$ is stationary if for any $k$, $n$ and $n_1, \dots, n_k$ in $\N$, the joint distribution of $(X_{n_1+n},\dots,X_{n_k+n})$ is equal to the joint distribution of $(X_{n_1},\dots,X_{n_k})$.
\end{defn}

\begin{defn}
 An $\cX$-valued stationary stochastic process $(X_n)_n$ is said to be ergodic if for every $\ell \geq 1$  and every pair of Borel sets $A,B \in \cX^{\ell}$,
\begin{align*}
 \lim_{n \to \infty} \frac{1}{n} \sum_{k=1}^{n} \PP\Bigl(  & (X_1,\dots,X_{\ell}) \in A, \, (X_{k+1},\dots,X_{k+\ell}) \in B \Bigr) \\
 & =  \PP\Bigl( (X_1,\dots,X_{\ell}) \in A \Bigr) \PP\Bigl((X_{1},\dots,X_{\ell}) \in B \Bigr).
\end{align*}
\end{defn}

\begin{defn}
A measurable dynamical system is a triple $(\cX,\cF,T)$, where $(\cX,\cF)$ is a measurable space and $T : \cX \to \cX$ is measurable. A topological dynamical system is a pair $(\cX,T)$, where $\cX$ is a topological space and $T : \cX \to \cX$ is a continuous map. In the study of topological dynamics, one often assumes that $\cX$ is compact and metrizable.
\end{defn}

\begin{defn}
 A measure-preserving system is a quadruple $(\cX,\cF,T,\mu)$, where $(\cX,\cF,\mu)$ is a measure space, $T: \cX \to \cX$ is measurable, and $\mu\bigl(T^{-1}(A)\bigr) = \mu(A)$ for each $A$ in $\cF$. In this case, we say that $T$ preserves the measure $\mu$ and $\mu$ is an invariant measure for $T$. For the purpose of this article, we will always assume that any invariant measure $\mu$ is a probability measure. Also, if $\cX$ is Polish and $\cF$ is the Borel $\sigma$-algebra, then we may refer to $(\cX,T,\mu)$ as a measure-preserving system.
\end{defn}

\begin{defn}
 A measure-preserving system $(\cX,\cF,T,\mu)$ is ergodic if $T^{-1}(A) = A$ implies $\mu(A) \in \{0,1\}$ for any $A$ in $\cF$. We may say that $T$ is ergodic for $\mu$, or we may say that $\mu$ is ergodic for $T$.
\end{defn}

 With the definitions given above, there is a correspondence between stationary stochastic processes and measure-preserving systems. Let us describe this correspondence as follows. Suppose $(X_n)_n$ is an $\cX$-valued stationary stochastic sequence, where $\cX$ is Polish. Let $\cY = \prod_n \cX$, equipped with the product $\sigma$-algebra induced by the Borel $\sigma$-algebra on $\cX$. Define $T : \cY \to \cY$ by the left shift: if $y = (x_n)_n$, then $(T(y))_n = x_{n+1}$. Kolmogorov's consistency theorem gives that there is a unique probability measure $\mu$ on $\cY$ with the same finite dimensional distributions as $(X_n)_n$. In this case, the stationarity of $(X_n)_n$ corresponds exactly to the invariance of $\mu$ with respect to $T$. Moreover, if $(X_n)_n$ is ergodic, then $\mu$ is ergodic for $T$.

 In the other direction, given any measure-preserving system $(\cX,T,\mu)$, we may define a stationary stochastic process as follows. For any Polish space $\cY$ and measurable map $f : \cX \to \cY$, let $X_n(\omega) = f(T^{n}(\omega))$. If $(\cX,T,\mu)$ is ergodic, then so is $(X_n)_n$.

Recall that an $\cX$-valued stochastic process $(X_n)_n$ is a Markov chain if for every $x$ in $\cX$, there exists a probability measure $\pi(x,\cdot)$ on $\cX$ such that for each measurable set $A$ in $\cX$, it holds that
\begin{equation*}
 \PP\Bigl( X_{n+1} \in A | X_1 = x_1, \dots, X_n = x_n \Bigr) = \pi(x_n,A).
\end{equation*}
In the model (\ref{generalobs})-(\ref{generaldyn}), if the dynamical noise process $(\delta_n)_n$ is assumed to be i.i.d., then both $(X_n)_n$ and $(X_n,Y_n)$ are Markov chains. This fact is particularly relevant in Sections \ref{SectionDynNoise} and \ref{SectionCombNoise}, where the process $(\delta_n)_n$ is assumed to be non-zero. Even in this case the process $(Y_n)_n$ may exhibit long-range dependencies. Setting the dynamical noise to zero in model  (\ref{generalobs})-(\ref{generaldyn}) can be thought of as as a very degenerate Markov chain, but it is not clear in this case how helpful the Markov perspective is, since even the process $(X_n)_n$ may exhibit long-range dependencies.

\subsection{Goals of statistical inference}

There are a variety of topics that can be considered part of ``statistical inference in dynamical systems." In the interest of providing context for this survey, let us mention the following topics:
\begin{enumerate}
  \item parameter estimation, model identification or reconstruction; 
  \item state estimation, filtering, smoothing, or denoising;
  \item feature estimation, where features often include invariant measures, dimensions, entropy, or Lyapunov exponents;
  \item prediction or forecasting;
  \item noise quantification, estimation, or detection.
\end{enumerate}

In this paper we focus almost exclusively on the problems of parameter inference, system identification or reconstruction. In the setting of (\ref{generalobs})-(\ref{generaldyn}), we pose the parameter estimation problem as follows. Suppose the family of dynamical systems can be parametrized by $T_a$, with parameter $a \in \cA$, as in \eqref{generaldyn}. Construct statistical procedures for estimating the parameter $a$, given observations $Y_1,Y_2,\cdots, Y_n$ from \eqref{generalobs}, and provide adequate theoretical support for the validity of the estimation procedure.

Of course, the boundaries between the problems mentioned above are often quite blurred. For example, if one can accurately estimate the hidden states $(X_k)_{k=0}^{n-1}$ from the data $(Y_k)_{k=0}^{n-1}$, then the problem of system identification often becomes significantly easier. For this reason, parameter inference methods often simultaneously attempt some version of state estimation or denoising.

\subsection{Organization of the paper}

We organize this survey according to which of the two types of noise in (\ref{generalobs})-(\ref{generaldyn}) are present (\textit{i.e.} non-zero). 
This organization is motivated by the observation that methods and results for parameter inference in dynamical systems tend to be specific 
to the type of noise assumed in the model.

The remainder of the paper is organized as follows. In Section \ref{SectionNoNoise} we describe some results relevant to inference for dynamical systems in the absence of noise. Section \ref{SectionObsNoise} contains a variety of proposed methods dealing with the case of dynamical systems contaminated by observational noise only. Section \ref{SectionDynNoise} deals with the case of only dynamical noise, and Section \ref{SectionCombNoise} addresses the setting of state space models, that is systems with both dynamical and observational noise. Lastly, we highlight some possibly interesting open questions in Section \ref{SectionOpenQuestions}.

 Ornstein and Weiss \cite{OW1991} have shown that in a certain sense it is impossible, in general, to tell the difference between observational and dynamical noise. In this sense, one might suggest that from the point of view of abstract ergodic theory, we should not make distinctions on the basis of the type of noise present. However, we are interested in finer properties than those captured by the equivalence relations considered in \cite{OW1991}, and therefore the distinction between observational and dynamical noise might still be useful for our purposes.

\subsection{Related surveys and books}

There have been many other reviews of topics related to the topics in this survey. An incomplete list of such reviews is the following: \cite{Berliner1992,BoylePetersen2011,ChatterjeeYilmaz1992,ER1985,Isham1993,Jensen1993,KostelichSchreiber1993,StemlerJudd2009,VTK2004}.
%Berliner (1992) \cite{Berliner1992}, Chatterjee and Yilmaz (1992) \cite{ChatterjeeYilmaz1992}, Eckmann and Ruelle (1985) \cite{ER1985}, Isham (1993) \cite{Isham1993}, Jensen (1993) \cite{Jensen1993}, Kostelich and Schreiber (1993) \cite{KostelichSchreiber1993}, Stemler and Judd (2009) \cite{StemlerJudd2009}, Voss, Timmer, Kurths (2004) \cite{VTK2004}.
Furthermore, let us mention the following books or monographs related to the topics in this survey: \cite{Abarbanel1996,BezruchkoSmirnov2010,ChanTong2001,DurbinKoopman2001,KantzSchreiber2004,Kifer1988,OW1991,Tong1990}.
%Abarbanel (1996) \cite{Abarbanel1996}, Bezruchko and Smirnov (2010) \cite{BezruchkoSmirnov2010}, Durbin and Koopman (2001) \cite{DurbinKoopman2001}, Kantz and Schreiber (2004) \cite{KantzSchreiber2004}, Kifer (1988) \cite{Kifer1988}, Tong (1990) \cite{Tong1990}, Ornstein and Weiss (1991) \cite{OW1991}.
The relevance of this survey is that we bring together approaches from many distinct fields and discuss them in a common statistical setting. In particular we discuss parameter estimation and inference  for the full range of noise settings. This perspective is rare since the different noise settings often correspond to different research areas such as deterministic dynamics or state space methods based on hidden Markov models. We bring these various approaches together and place them in a common context. Inference in dynamical systems for a variety of contexts was discussed in  Berliner (1992) \cite{Berliner1992}, and our survey can be thought of as an updated and greatly expanded version of this work.

%%%%%%%%%%%%%%%%%%%%%%%%%%%%%%%%%%%%%%%%%%%%%%%%%%%%%%%%%%%%%%%%%%%%
%%%%%%%%%%%%%%%%%%%%%%%%%%%%%%%%%%%%%%%%%%%%%%%%%%%%%%%%%%%%%%%%%%%%
\section{No noise} \label{SectionNoNoise}

If no noise is present in the model (\ref{generalobs})-(\ref{generaldyn}), then we have the following situation:
\begin{align}
Y_n & = f_a(X_n) \label{nonoiseobs} \\
X_{n+1} & = T_a(X_n), \label{nonoisedyn}
\end{align}
where $T : \cA \times \cX \to \cX$ is a parametrized family of maps and $f : \cA \times \cX \to \cY$ is a parametrized family of observation functions. For a fixed parameter value $a$, the model (\ref{nonoiseobs})-(\ref{nonoisedyn}) is one of the classical objects of study in dynamical systems and ergodic theory (for general references on dynamical systems and ergodic theory, see  \cite{BrinStuck2002,KatokHasselblatt1995,Petersen1989,Walters1982}).

\subsection{Non-parametric system reconstruction from direct observations} \label{SectionNoNoiseReconstructionDirectObservations}

Here we consider non-parametric estimation of a map $T$ from direct observation of a single trajectory. Although the methods discussed in this section do not directly involve parameter estimation, they are nonetheless relevant for parameter estimation, since any non-parametric method for estimation of a map immediately yields a method of parameter estimation if the map to be estimated comes from a parameterized family.

Let us first consider a case when the system can be successfully reconstructed from observations. If $\cX$ is a manifold, $T$ is continuous, the trajectory $(x_n)_n$ is dense in $\cX$, and we observe the trajectory directly (\textit{i.e.} the observations $(y_n)_n$ satisfy $x_n = y_n$), then $T$ can be consistently estimated from $(y_n)_n$ using locally linear functions of the data. More precisely, let us state a result from \cite{AdamsNobel2001} justifying this statement in the case $\cX = [0,1]$. Let $\lambda$ be Lebesgue measure on $[0,1]$. The map $T : [0,1] \to [0,1]$ is said to be an $E\{I_j,\alpha_j\}$-map if there exists at most countably many disjoint open intervals $I_j$ and real numbers $\alpha_j$ such that $\lambda( \cup I_j ) = 1$ and $f'(x) = \alpha_j$ for all $x$ in $I_j$.

\begin{prop}[\cite{AdamsNobel2001}]
Let $T$ be an $E\{I_j,\alpha_j\}$-map. Suppose the observed trajectory $(x_n)_n$ is dense in $[0,1]$. Then there exists a sequence of estimates $\hat{T}_n$ of $T$ such that for almost every $x$ in $[0,1]$, it holds that $\hat{T}_n(x) = T(x)$ for all but finitely many $n$. In particular, $\hat{T}_n$ converges to $T$ pointwise almost everywhere, and $\lambda( \{ x : \hat{T}_n \neq T(x) \} )$ tends to zero.
\end{prop}

To get an idea about how to prove this proposition, notice that for any two consecutive points $x_n$ and $x_{n+1}$ in the trajectory, the pair $(x_n,x_{n+1})$ lies on the graph of $T$. Therefore one may estimate $T$ by linearly interpolating between neighboring points on the graph of $T$.

When the map $T$ is not assumed to be continuous but only measurable, estimation of $T$ from discrete observations of a single trajectory has been carried out by Adams and Nobel \cite{AdamsNobel2001}. In this work, the map $T$ is assumed to preserve a Borel probability measure $\mu$ on $\cX$, and the system $(X,\mu,T)$ is assumed to be ergodic. Their main result may be stated as follows. 

\begin{thm}[\cite{AdamsNobel2001}] \label{AdamsNobelThm} Let $\mu_0$ be a reference probability measure on $\cX$ that is assumed to be ``known.'' Also assume that there is a ``known" constant $M$ such that $1/M \leq d\mu/d\mu_0 \leq M$. Let $\Meas(\cX)$ denote the space of measurable functions from $\cX$ to $\cX$. Then there is a estimation scheme $(T_n)_n$ (whose definition uses $M$ and $\mu_0$), where $T_n : \cX^n \to \Meas(\cX)$, such that for $\mu_0$-a.e. initial condition $x_0$, the map $T_n(x_0,\dots,x_{n-1})$ converges to $T$ in a weak topology (\textit{i.e.} $\mu\bigl(T_n^{-1}(A) \bigtriangleup T^{-1}(A)\bigr)$ tends to zero as $n$ tends to infinity for each Borel set $A$).
\end{thm}

The estimation scheme $(T_n)_n$ that appears in \cite{AdamsNobel2001} is constructed using an adaptive histogram method, which we discuss below. This paper also shows that under the same hypotheses the conclusion of the theorem is false if one requires that $\mu( \{x \in \cX : T_n(x) \neq T(x)\})$ tends to zero as $n$ tends to infinity.

Here we give an idea of the estimation scheme used in the proof of Theorem \ref{AdamsNobelThm}. The histogram method described here is actually from \cite{Nobel2001}, which is very similar in spirit to the method used in the proof of Theorem \ref{AdamsNobelThm}. Assume that $\cX \in \R^d$, and we fix a refining sequence $(\pi_k)_k$ of finite partitions of $\cX$ with some additional properties (see \cite{AdamsNobel2001} for details). Let $\pi_k(x)$ denote the cell in $\pi_k$ containing $x$. Given the first $n$ terms of the trajectory $(x_j)_{j=0}^{n-1}$, let
\begin{equation*}
\phi_{n,k}(x) = \frac{ \sum_{j=0}^{n-1} x_{j+1} I_{\{x_j \in \pi_k(x) \} } }{ \sum_{j=0}^{n-1} I_{\{x_j \in \pi_k(x)\} }},
\end{equation*}
where $I_{\{x_j \in \pi_k(x)\}}$ is the indicator function of the event that $x_j$ is in $\pi_k(x)$, and if the cell $\pi_k(x)$ contains no points $x_j$, then $\phi_{n,k}(x) = 0$. Now consider the empirical loss of $\phi_{n,k}$:
\begin{equation*}
\Delta_{n,k} = \biggl( \frac{1}{n} \sum_{j=0}^{n-2} (\phi_{n,k}(x_j) - x_{j+1})^2 \biggr)^{1/2}.
\end{equation*}
The estimates $\hat{T}_n$ of $T$ are adaptively chosen from among the $\phi_{n,k}$ according to $\Delta_{n,k}$ (using $\mu_0$ and $M$). This method has the advantage that it works in quite a general setting (the only assumptions involve ergodicity and the Radon-Nikodym derivative with respect to a reference measure). On the other hand, it relies on the ergodic theorem for convergence, and therefore it appears very unlikely that it would have any general speed of convergence.

\subsection{Non-parametric system reconstruction from general observations}

In this section we consider approaches to system reconstruction when the observations $(y_n)_n$ are not necessarily equal to the trajectory $(x_n)_n$. There is a vast amount of literature on the technique of system reconstruction via delay coordinate embeddings. These system reconstructions may be thought of as non-parametric inference of dynamical systems. 
%In this setting one is given some observations of a trajectory and the goal is to identify the structure of the underlying system. 
Delay coordinate embeddings are a well-studied inference procedure to reconstruct dynamical systems that satisfy certain conditions. In this section we define delay coordinate embeddings, mention some of the main uses of these techniques, and provide some representative theorems that provide conditions under which these methods work.

The eventual goal of delay coordinate embedding techniques is typically feature estimation, which we summarize as follows. If the underlying map $T$ and the observation function are both smooth, then under generic conditions, a delay coordinate embedding allows one to construct a smooth map $\widetilde{T}$ such that $\widetilde{T}$ is related to $T$ by a smooth change of coordinates. Under this scenario, $T$ and $\widetilde{T}$ will share many features, including entropy, Lyapunov exponents, and fractal dimensions of corresponding invariant measures. As these features are considered important in many physical settings, such delay coordinate reconstructions have been extensively studied.

To be specific, we consider a smooth map $T : \cX \to \cX$ of a manifold $\cX$, with a smooth observation function $f : \cX \to \R$. The data are assumed to be generated as follows: there is a trajectory $(x_n)_n$ such that $x_{n+1} = T(x_n)$, and we observe the data $(y_n)_n$ such that $y_n = f(x_n)$. The original idea to use delay coordinate embeddings to construct a system equivalent to $(\cX,T)$ from the observations is due to Ruelle, at least according to the influential paper \cite{PCFS1980}. 
\begin{defn}
A delay coordinate mapping of $\cX$ into $\R^m$ is a mapping $F : \cX \to \R^m$ such that
\begin{equation*}
F(x) = (f(x),f \circ T^\tau(x), \dots, f \circ T^{\tau(m-1)}(x) ),
\end{equation*}
for some natural number $\tau$. The mapping $F$ is said to be an embedding if it is a diffeomorphism from $\cX$ to its image $F(\cX)$, that is if $F$ is a smooth injection and has a smooth inverse. 
\end{defn}
The well-known theorem of Takens \cite{Takens1981} (often called the Takens Embedding Theorem) may be stated as follows. 
\begin{thm}[\cite{Takens1981}]
If $T$, $f$, and $\tau$ satisfy certain genericity conditions and $m > 2 \dim(\cX)$, then $F$ is an embedding.
\end{thm}
Let $\tilde{\cX} = F(\cX)$ and $\tilde{T} = F \circ T \circ F^{-1}$. The fact that $F$ is an embedding means that the system $(\cX,T)$ is related to the system $(\tilde{\cX},\tilde{T})$ by a smooth change of coordinates (given by $F$). In particular, invariants of $(\cX,T)$ that depend on the differential structure of $T$ (such as Lyapunov exponents or fractal dimensions of attractors) are equal to those of the system $(\tilde{\cX},\tilde{T})$.

In particular, given the data $(y_k)_{k=0}^{n-1}$, we may build time series data $(s_k)_{k=0}^{n-1-\tau(m-1)}$ for the system $(\tilde{\cX},\tilde{T})$ as follows: for $k = 0, \dots, n-1-\tau(m-1)$, let
\begin{equation*}
 s_k = (y_k,y_{k+\tau},\dots,y_{k+\tau(m-1)}).
\end{equation*}
Then the new time series $(s_k)_k$ may be used to estimate invariant features of $(\tilde{\cX},\tilde{T})$, which will be the same as those features of $(\cX,Q)$.

Takens's theorem has been generalized in various directions, such as filtered delay embeddings (see \cite{SYC1991}, for example) or delay embeddings for stochastic systems (see \cite{SBDH2003}), but we do not attempt to record all such results. However, the following generalization, due to Sauer, Yorke, and Casdagli, bears mentioning.
\begin{thm}[\cite{SYC1991}]
 Let $A$ be a compact subset of $\cX$ with box-counting dimension $d$. Let $m > 2d$. Suppose $T$, $f$, $\tau$, and $A$ satisfy certain genericity conditions. Then the delay coordinate map $F$ given above is an injection on $A$ and an immersion on each compact subset of any smooth manifold contained in $A$. 
\end{thm}
The advantage of this theorem over the Takens theorem is that the relevant dimension $d$ might be less than the ambient dimension of $\cX$, in which case the number of coordinates $m$ required in the embedding space may be less than the number of coordinates required by Takens's theorem.

In order to use the delay coordinate method given only the data $(y_k)_{k=0}^{n-1}$, one must choose an appropriate dimension $m$ and an appropriate lag $\tau$. A variety of statistical techniques have been proposed to estimate the dimension $m$ and find a suitable lag $\tau$ (for example, see the book \cite{KantzSchreiber2004} or the collection \cite{MeesCollection2001}), but further pursuit of these topics lies outside the scope of this survey.

\subsection{Results from ergodic theory} \label{SectionErgTheory}
In this section, we state some results  from ergodic theory that are relevant for parameter inference. %Mention here: \cite{GH2008,OW1990,OW1991,OW2007,RyabkoRyabko2010}, along with \cite{HochmanQuas}. 

One of the most general results in this area is due to Ornstein and Weiss \cite{OW1990}. In this work, the authors consider the problem of estimation of stationary ergodic processes. (Note that in the setting of (\ref{nonoiseobs})-(\ref{nonoisedyn}), if $X_0$ is distributed according to an ergodic invariant measure for $T_a$, then the observation process $(Y_n)_n$ satisfies exactly these conditions.) To make this problem precise, they consider the $\overline{d}$ metric on the space of such processes. Their main results may be stated as follows. First, they construct a procedure which, given a realization $(X_k)_{k=0}^{n-1}$ of a process $(X_k)_k$ constructs a process $Z^n = (Z^n_k)_k$. Then they show that the sequence of processes $(Z^n)_n$ converges to $(X_k)_k$ in the $\overline{d}$ metric if and only if $(X_k)_k$ is Bernoulli. Thus, they have shown that there is a consistent estimation procedure for the class of Bernoulli processes. Furthermore, they show that no estimation procedure can be consistent for the class of all stationary ergodic processes. 

In another direction, Ornstein and Weiss \cite{OW2007} show that entropy is the only finitely observable invariant in the following sense. Let $J$ be a function from the class of finite-valued stationary ergodic processes to a complete separable metric space such that $J$ is constant on isomorphism classes. The main result of \cite{OW2007} states that if $J$ is finitely observable, then it must be a continuous function of the entropy. This result shows that there are strong restrictions on the possibilities for inference of isomorphism invariants.

Gutman and Hochman \cite{GH2008} extend the results in \cite{OW2007} in several ways. They give several rich families of classes $\mathcal{C}$ of stationary ergodic processes such that if $J$ is a finitely observable invariant on $\mathcal{C}$, then $J$ is constant. They also show that for every finitely observable invariant $J$ on the class of irrational circle rotations, $J$ is constant on the processes arising from a full measure set of angles. In particular, there is no finitely observable invariant for irrational rotations which is complete.

%Hochman and Quas \cite{HochmanQuas} consider the class of circle rotations, in which a map is parametrized by its angle of rotation $a \in [0,1/2)$. Let $\cX = S^1$, where the map $Q_a: S^1 \to S^1$ is given by $Q_a(x) = x+a \mod 1$. They prove two results, one of which is positive and the other negative. The negative result states that there is no estimation scheme that will provide consistent estimation of $a$ for all continuous observation functions. On the other hand, the positive result constructs an estimator for $a$ that is consistent for any $C^{1/2+\beta}$ observation function, for any $\beta > 0$. Together, these results highlight the subtle role that the regularity of the observation function plays in ergodic theory and dynamical systems.

There is a large body of work, often categorized as smooth ergodic theory, that seeks to understand the statistical properties of smooth (or piecewise smooth) dynamical systems. The typical setting is that one has a compact Riemannian manifold $M$ and a smooth self-map $f : M \to M$. The manifold typically has a distinguished probability measure $\lambda$, which one may think of as volume measure on the manifold. The goal is to understand the asymptotic behavior of the trajectory $\{f^n(x)\}_n$ for $\lambda$-a.e. $x$. For a wide class of such systems \cite{Young1998}, often called (non-uniformly) hyperbolic systems, there is an invariant measure $\mu$ on $M$ such that for $x$ in a set of positive $\lambda$-measure, the trajectories in $x$ equidistribute with respect to $\mu$. In such cases, the measure $\mu$ is said to be a \textit{physical} measure. Often the measure $\mu$ has some additional properties (it has no zero Lyapunov exponents and absolutely continuous conditional measures with respect to $\lambda$ on unstable manifolds), and in this case $\mu$ may be called an SRB (Sinai-Ruelle-Bowen) measure \cite{Young2002}. The ergodic theory of SRB measures is fairly well-studied, and many of their statistical properties have been analyzed.

A related topic that has seen a great deal of attention recently is concentration inequalities for dynamical systems  
\cite{Chazottes2012,CCRV2009,CCS2005_1,CCS2005_2,ChazottesGouezel2011}. These inequalities are used to study
 the fluctuations of observables for dynamical systems and have been shown to hold for sufficiently regular observables and 
 a wide class of non-uniformly hyperbolic dynamical systems. Using these inequalities, it is possible to perform statistical estimation 
 of various features of the dynamical system. See the survey \cite{Chazottes2012} for more details and precise statements.

\subsection{Parameter inference via synchronization and control}

Synchronization-based approaches to parameter estimation have appeared quite often in the physics and control systems literature \cite{ACJ2008,MaybhateAmritkar1999,Parlitz1996,QBCKA2009,YCCLP2007}. In situations when these methods are used, it is common that no particular noise model is assumed. Indeed synchronization-based approaches are typically described as parameter inference methods in the noiseless setting, although they may be applied in other settings. The main idea of synchronization-based methods is to insert a ``control'' term in the defining equations of the system that allows one to incorporate the data. The parameter estimation may then be framed as a large optimization procedure in which one tries to find trajectories of the system which are close to the data.

The topic of parameter estimation in a noiseless setting is discussed directly in the work of Abarbanel, Creveling, Farsian, and Kostuk \cite{ACFK2009}, and we review their approach in this section. The main issue in this context is that one only has access to the observations $(Y_n)_n$, which might ``hide" some information about the underlying system. The approach taken in \cite{ACFK2009} involves synchronization of the observations and the output of a model over the relevant time window. This approach may be summarized as follows.

Suppose that $\cX$ is in $\R^{d}$ and the system (\ref{nonoiseobs})-(\ref{nonoisedyn}) has the following form:
\begin{align*}
Y_n &= X_{n,1} \\
X_{n+1,i} & = T_{a,i}(X_n),
\end{align*}
where $X_{n,i}$ denotes the $i$-th coordinate of $X_n$. The synchronization approach taken in \cite{ACFK2009} is to add a ``control" term of the form $k(Y_n - X_{n,1})$ to first coordinate of the model as follows:
\begin{align*}
\tilde{X}_{n+1,1} & = T_{a,1}(\tilde{X}_n) + k (Y_n - \tilde{X}_{n,1}) \\
\tilde{X}_{n+1,i} & = T_{a,i}(\tilde{X}_n), \quad i > 1.
\end{align*}
For $k >0$ large enough, the data $Y_n$ and the first coordinate $\tilde{X}_{n,1}$ of the model trajectory will ``synchronize." With a fixed $k$, the authors propose to estimate the parameter $a$ and the initial state $X_0$ by minimizing the following function:
\begin{equation*}
C(a,X_0) = \sum_{j=0}^{n-1} (Y_n - \tilde{X}_{n,1})^2,
\end{equation*}
where the trajectory $\tilde{X}_n$ is computed starting at $\tilde{X}_0 = X_0$. The purpose of adding the control term is to regularize the function $C$ so that its minimum may be found efficiently. Of course, the trajectory $\tilde{X}_n$ associated with this minimum is not a true trajectory of the original system. Therefore the authors propose a synchronization method that allows the parameter $k$ to depend on time. In other words, they propose to minimize the cost function
\begin{equation*}
C(a,X_0) = \sum_{j=0}^{n-1} (Y_j - \tilde{X}_{j,1})^2 + k_j^2,
\end{equation*}
subject to the constraints
\begin{align*}
\tilde{X}_{j+1,1} & = T_{a,1}(\tilde{X}_j) + k_j (Y_j - \tilde{X}_{j,1}) \\
\tilde{X}_{j+1,i} & = T_{a,i}(\tilde{X}_j), \quad i > 1.
\end{align*}
Although this method has been observed to work sufficiently well in practice \cite{ACFK2009}, we remark that to the best of our knowledge there are no theoretical guarantees regarding the consistency or performance of this method.

%%%%%%%%%%%%%%%%%%%%%%%%%%%%%%%%%%%%%%%%%%%%%%%%%%%%%%%%%%%%%%%%%%%%
%%%%%%%%%%%%%%%%%%%%%%%%%%%%%%%%%%%%%%%%%%%%%%%%%%%%%%%%%%%%%%%%%%%%
\section{Observational noise only} \label{SectionObsNoise}

If only observational noise is present in the model (\ref{generalobs})-(\ref{generaldyn}), then the system (\ref{generalobs})-(\ref{generaldyn}) reduces to the following situation:
\begin{align}
Y_n & = f_a(X_n,\epsilon_{n}) \label{obsnoiseobs} \\
X_{n+1} & = T_a(X_n), \label{obsnoisedyn}
\end{align}
where $(\epsilon_n)_n$ is a noise process, $T : \cA \times \cX \to \cX$ is a parametrized family of maps, and $f : \cA \times \cX \times \cN \to \cY$ is a parametrized family of noisy observation functions. Multiple authors explicitly argue for consideration of the observational noise model. For example, Judd \cite{Judd2003_2} states that ``the reality is that many physical systems are indistinguishable from deterministic systems, there is no apparent small dynamic noise, and what is often attributed as such is in fact model error." Furthermore, Lalley and Nobel \cite{LalleyNobel2006} remark that ``estimation in the observational noise model has not been broadly addressed by statisticians, though the model captures important features of many experimental situations."

A distinguishing feature of the observational noise model is that the process $(X_n)_n$ is deterministic, and therefore in general it exhibits a long-range dependence structure. Furthermore, this long-range dependence is still present beneath the noise in the observation process $(Y_n)_n$. Such dependencies imply that traditional statistical estimation techniques do not apply and may not work. As Lalley and Nobel state in \cite{LalleyNobel2006}, ``though some features of denoising can be found in more traditional statistical problems such as errors in variables regression, deconvolution, and measurement error modeling (c.f. \cite{CRS1995}), other features distinguish it from these problems and require new methods of analysis.'' In particular, they cite the facts that the covariates $X_n$ are deterministically related (as opposed to i.i.d. or mixing), the noise is often bounded (as opposed to Gaussian), and the noise distribution itself is often unknown.

%\textbf{Standing assumption in this Section.} For the remainder of this section, we assume that $\cY = \cX \subset \R^d$ and $f_a(x,\epsilon) = x + g_a(\epsilon)$, where $g_a$ is a function specifying the effect of the additive observational noise.

\begin{example} \label{ExampleLogistic}
Let $\cX = [0,1]$, and let $T_a : \cX \to \cX$ be given by $T_a(x) = ax(1-x)$, with $a$ in $\cA = [0,4]$. This family of maps, known as the \textit{logistic family}, has been extensively studied in a variety of settings. For $a \in [0,1]$, it is known that for all $x$ in $[0,1]$, the iterates $T_a^n(x)$ tend to $0$ as $n$ tends to infinity. We say that a parameter value $a$ has an attracting periodic orbit $\{p_0, \dots, p_{N-1}\}$ if $T_a(p_i) = p_{i+1}$ (with indices interpreted modulo $N$) and $|(T_a^{N})'(p_i)| < 1$. For such parameter values, the iterates $T_a^n(x_0)$ of Lebesgue almost every initial point $x_0$ will tend to the periodic orbit $\{p_0,\dots, p_{N-1}\}$ as $n$ tends to infinity. It is known \cite{GS1997, Lyubich1994} that the set of parameter values that have an attracting periodic orbit is open and dense in $[0,4]$. On there other hand, there are parameter values that give rise to very different asymptotic dynamics. In particular, we say that a parameter value $a$ has an absolutely continuous invariant measure (acim) $\mu_a$ if $\mu_a$ is absolutely continuous with respect to Lebesgue and $\mu_a$ is an invariant measure for $T_a$. In such cases, it can be shown that the iterates $T^n_a(x_0)$ of Lebesgue almost every initial point $x_0$ equidistribute with respect to $\mu_a$. Intuitively, the presence of $\mu_a$ produces seemingly stochastic behavior, which is often referred to as chaos.  Jakobson showed in \cite{Jakobson1981} that the set of parameter values that have an acim has positive measure in $[0,4]$, and Lyubich eventually showed in \cite{Lyubich2002} that Lebesgue almost every parameter in $[0,4]$ either has an attracting periodic orbit or an acim.

In most of the papers cited in this section, this family of maps is taken as a standard testing ground for parameter estimation methods. Generally, it is assumed that the observational noise is additive (\textit{i.e.} $f_a(x,\epsilon) = x + \sigma(a)\epsilon$).
\end{example}

\subsection{Noise reduction}
One basic approach to parameter estimation in the observational noise case is to reduce the noise and then apply parameter estimation methods. If the noise can be uniformly and sufficiently reduced, then these approaches will be approximately as successful as the estimation method applied to the noiseless case. For example, the positive results in \cite{Lalley1999,Lalley2001,LalleyNobel2006} might be combined with a parameter estimation method in order to produce consistent estimates. Among the results contained in these works, the main positive result of \cite{LalleyNobel2006} is the most general, and we state it as follows. 

A homeomorphism $F$ of a compact metric space $(\Lambda,d)$ is said to be expansive with separation threshold $\Delta$ if for every $x \neq y$ in $\Lambda$, there exists $n$ in $\Z$ such that $d(F^n(x),F^n(y)) >\Delta$. In the work \cite{LalleyNobel2006}, the authors consider an initial condition $x$ and let $x_i = F^i(x)$. Also, they define a particular denoising algorithm which, given noisy additive noisy observations $(y_i)_{i=0}^{n-1}$, produces estimates $\hat{x}_{i,n}$ of the true states $x_i$. In this context, the main positive result may be stated as the following theorem.
\begin{thm}[\cite{LalleyNobel2006}] \label{ThmLalleyNobelPositive}
 Let $F : \Lambda \to \Lambda$ be an expansive homeomorphism with separation threshold $\Delta > 0$. Suppose that the noise process $(\epsilon_n)_n$ satisfies $|\epsilon_n| \leq \Delta/5$ for every $n$. If $k = k(n) \to \infty$ and $k / \log(n) \to 0$ as $n$ tends to infinity, then
\begin{equation*}
 \frac{1}{n-2k} \sum_{i=k}^{n-k} |\hat{x}_{i,n} - x_i| \to 0, \quad \text{ as } n \to \infty
\end{equation*}
with probability $1$ for every initial point $x$ in $\Lambda$ (with respect to any $F$ invariant Borel probability measure).
\end{thm}
By allowing a slight modification to their estimation scheme, the authors also show that under the same hypotheses
\begin{equation*}
 \max_{\log(n) \leq i \leq n - \log(n)} |\hat{x}_{i,n} - x_i| \to 0, \quad \text{ as } n \to \infty
\end{equation*}
with probability $1$ for almost every initial point $x$ in $\Lambda$.

Of course, the task of removing the noise might itself be difficult or in some cases even impossible, as witnessed by the negative results in \cite{Lalley1999,Lalley2001,LalleyNobel2006} and the related results in \cite{Judd2003,Judd2007,JuddSmith2001,Judd2004}. Here we state the main negative result in \cite{LalleyNobel2006}.  A pair of points $x$ and $x'$ is said to be strongly homoclinic for the homeomorphism $F$ if
\begin{equation*}
\sum_{k \in \Z} d(F^k(x),F^k(x')) < \infty.
\end{equation*}

\begin{thm}[\cite{LalleyNobel2006}] \label{ThmLalleyNobelNegative}
 Suppose the stationary distribution for the noise process $(\epsilon_n)_n$ is unbounded (or has sufficiently large support). If $x$ and $x'$ are strongly homoclinic, then for every measurable function $\phi : \prod_n \cX \to \R^d $, 
\begin{equation*}
\EE\biggl[| \phi((y_n)_n) - x| - | \phi((y'_n)_n) - x'|\biggr] > 0.
\end{equation*}
\end{thm}
In other words, even with access to the entire observation sequence, any state estimation or denoising scheme will fail with positive probability.

Lalley and Nobel also point out that despite the fact that in such cases the problem of asymptotic denoising is impossible, it might still be possible to obtain consistent parameter estimates. In fact, they state that such examples might provide an interesting avenue for further study (see Question \ref{QuestionLalleyNobel}).

In addition to the works mentioned so far in this section, the following works discuss the problem of denoising or smoothing data in the presence of only observational noise: \cite{Davies1994,GSHT2010,KostelichSchreiber1993,KostelichYorke1990,MGA2009,MGA2010,Sauer1992}.

\subsection{Introduction to likelihoods and related methods} \label{SectionIntroToLikelihoods}

We begin with the work of Berliner \cite{Berliner1991,Berliner1992}, sets the stage for most of the work that has followed. %(The work of MacEachern and Berliner \cite{MacEachernBerliner1995} deals with a non-stationary situation in which the parameters change over time, which lies outside the scope of this survey.)
In these works, the author is mostly concerned with the observational noise setting (\ref{obsnoiseobs})-(\ref{obsnoisedyn}). The likelihood function is given by
\begin{equation*}
L(x_0,a) = p(y_0^{n-1} | x_0,a),
\end{equation*} 
where $p(y_0^{n-1}|x_0,a)$ denotes the likelihood of observing $y_0^{n-1}$ given the parameter choice $a$ and the true initial condition $x_0$ (\textit{i.e.} $p( \cdot | x_0,a)$ is the probability density for the observation process conditional on $x_0$ and $a$). The maximum likelihood (ML) method for estimating the parameter $a$ amounts to defining the following maximum likelihood estimator (MLE):
\begin{equation} \label{EquationMLE}
\hat{a}_n = \underset{a}{\argmax} \max_{x_0} L(x_0,a).
\end{equation}
It will be useful to find an explicit form for the likelihood function in the case that (I) the observational noise sequence $(\epsilon_n)_n$ is assumed to be i.i.d. normal with zero mean and unit variance, and (II) the observation function $f_a$ takes the form $f_a(x,\epsilon) = x + \sigma(a) \epsilon$. The function $\sigma(a)$ allows one to set the variance of the noise according to the parameter $a$. In this case, we have
\begin{equation*}
L(x_0,a) = \Bigl(\sigma(a) \sqrt{2\pi} \Bigr)^{-n} \exp\Biggl( -  \sum_{k=0}^{n-1} (y_k - T_a^k(x_0))^2/(2\sigma^2(a))  \Biggr)\,
\end{equation*}
and the corresponding log-likelihood function is given by
\begin{equation} \label{LLik}
\log L(x_0,a) = -n \log \Bigl(\sigma(a) \sqrt{2\pi}\Bigr) - \sum_{k=0}^{n-1} (y_k - T_a^k(x_0))^2/(2\sigma^2(a)).
\end{equation}
A significant portion of the work on parameter estimation following Berliner has involved optimization of this log likelihood function, even when the noise is not necessarily Gaussian and thus its interpretation as a log likelihood function is no longer valid. 

As discussed in \cite{PS2004}, no existing statistical results apply to the ML method in this setting. With the above notation, the main difficulty in the current setting is that $T_a^k$ is a non-stationary function of $k$. Standard statistical results on the performance of the ML method apply when the likelihood function has no such dependence on $k$ (or is periodic with respect to $k$), but these results do not apply \textit{a priori} in the current setting.

The Bayesian approach  assumes a prior distribution (density) for $x_0$ and $a$, written as $\pi(x_0,a)$. Given the data $y_0^{n-1}$, the posterior distribution is then
\begin{equation*}
\pi(x_0,a|y_0^{n-1}) = \frac{p(y_0^{n-1}|x_0,a) \pi(x_0,a)}{ \int p(y_0^{n-1} | x,a) \pi(x,a) \, dx da}.
\end{equation*}

In these basic definitions, Berliner considers three main methods of parameter estimation: maximum likelihood estimation, minimization of a cost function (which is often chosen to be the negative of the log likelihood function) and Bayesian estimation. 
%All of these methods involve consideration of the likelihood function $L(x_0,a)$. 
One of Berliner's main points is that when the system (\ref{obsnoiseobs})-(\ref{obsnoisedyn}) is chaotic, the likelihood function will also typically be chaotic, in the sense that it will be extremely jagged.  The rough nature of these likelihood functions makes all three of the above methods of statistical estimation computationally very expensive, and much of the work following Berliner has been motivated by the need to mitigate this difficulty. Beyond these computational difficulties, we would like to emphasize that to our knowledge there are no general results concerning the consistency of any of these likelihood-based methods. 

%\begin{figure}[h!]
%  \caption{Characteristics for Part I (2) c}
%  \centering
%    \includegraphics[scale = .8]{prob4c_p37part2.pdf}  \label{Char_I2c}
%\end{figure}

\subsection{Variations on likelihood based methods} \label{SubSectionVariations}

A common method of parameter estimation in practice is to minimize some cost function $C$ with respect to the parameters. Given the observations $(y_k)_{k=0}^{n-1}$, such methods employ the following estimators:
\begin{equation*}
 \hat{a}_n = \argmin_a \min_{x_0} C\left(x_0,a,(y_k)_{k=0}^{n-1}\right),
\end{equation*}
where $C(x_0,a,(y_k)_{k=0}^{n-1})$ somehow measures the discrepancy of the observations and the system trajectory having  parameter $a$ and initial state $x_0$.

As we mentioned in the previous section, the most basic cost function is the least squares cost function
\begin{equation} \label{LScost}
C_{_{LS}}\left(x_0,a,(y_k)_{k=0}^{n-1}\right) = \sum_{k=0}^{n-1} (y_{k} - T_a^k(x_0) )^2.
\end{equation}

\begin{figure}[ht]
  \centering
    \includegraphics[scale = .8]{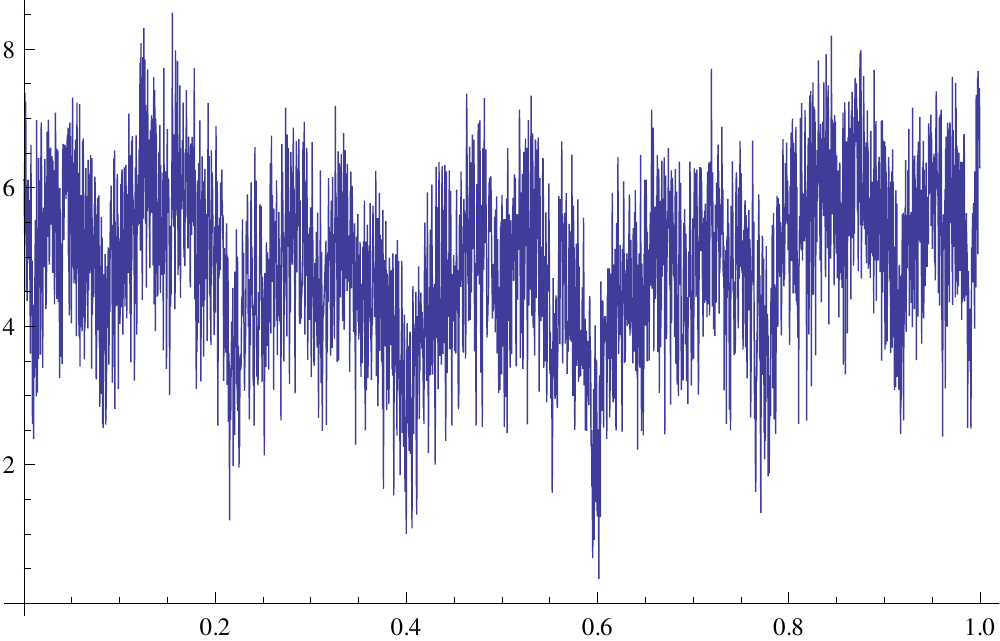}  \label{LSCostFunctionx}
 \caption{Least Squares cost function for $x_0$ in logistic family as a function of $x \in [0,1]$ given $n=20$ observations, true initial value $x_0 =.4$ and true parameter $a = 4$.}
\end{figure}

\begin{figure}[ht]
  \centering
    \includegraphics[scale = .8]{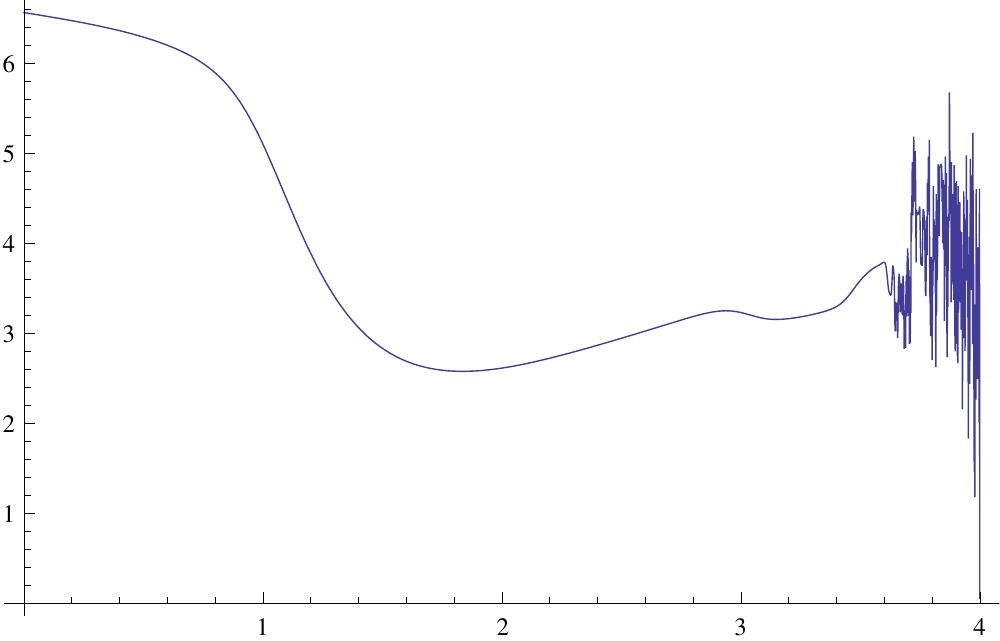}  \label{LSCostFunctiona}
 \caption{Least Squares cost function for parameter $a$ in logistic family as a function of $a \in [0,4]$ given $n=20$ observations, true initial value $x_0 =.4$ and true parameter $a = 4$.}
\end{figure}

Perhaps due to the sensitive dependence of $C_{_{LS}}$ on $x_0$ and the additional computational expense incurred by minimizing $C_{_{LS}}$ over $x_0$, several authors considered minimization of a one-step least squares cost function, given by
\begin{equation} \label{OSLScost}
C_{_{OSLS}}(a,(y_k)_{k=0}^{n-1}) = \sum_{k=0}^{n-2} (y_{k+1} - T_a(y_k))^2,
\end{equation}
which does not depend on any initial condition $x_0$. This cost function may appear to be the familiar least squares function from regression analysis, but as Kostelich \cite{Kostelich1992} recognized, it suffers from the problem of errors in variables (\textit{c.f.} \cite{CRSC2006,Fuller2006}). The problem of errors in variables is 
that the errors are not independent as is assumed by the cost function. Viewing $C_{_{OSLS}}$ from the perspective of traditional regression, we see that $y_k$ appears to play the role of the independent variable and $y_{k+1}$ plays the role of the dependent variable, but both $y_k$ and $y_{k+1}$ contain noise according to the model (\ref{obsnoiseobs})-(\ref{obsnoisedyn}). It is well-known that the problem of errors in variables can lead to asymptotically biased results, and therefore we should not expect minimization of $C_{_{OSLS}}$ to give consistent estimates of the parameter $a$.

In response to the errors in variables problem, Jaeger and Kantz \cite{JaegerKantz1996,KantzJaeger1997} propose a ``solution" of the problem, which amounts to minimizing the following cost function that has since gone by the name ``total least squares'' cost function:
\begin{equation} \label{TLScost}
C_{_{TLS}}(a,(y_k)_{k=0}^{n-1})  = \sum_{k=0}^{n-1} \min_{y \in \cX} |(y_k,y_{k+1}) - (y,T_a(y))|^2.
\end{equation}
Note that this approach essentially ignores the dynamics altogether, and instead focuses on minimizing the sum of orthogonal distances between the graph of $T_a$ and the points $(y_k,y_{k+1})$ in $\cX \times \cX$. In order to include some aspect of the dynamics, they further modify their cost function to find local shadowing trajectories by considering cost functions of the form
\begin{equation} \label{modTLScost}
C_{_{MTLS}}(a) = \sum_{k=0}^{n-s-1} \min_{y} |(y_k,\dots,y_{k+s}) - (y,\dots,T_a^{s}(y))|.
\end{equation}
Here $s$ is a parameter of the method; it is the number of steps over which one considers the local shadowing trajectories. If one asks for global shadowing trajectories, corresponding to $s = n-2$, then this modified total least squares cost function is equivalent to the original least squares cost function $C_{_{LS}}$.

McSharry and Smith \cite{McSharrySmith1999} consider the one step cost function $C_{_{OSLS}}$ given by (\ref{OSLScost}). They prove that in the case of the logistic map with a specific parameter value, the minimization of this cost function produces biased estimates, even with infinitely many observations. Their proposed solution involves minimizing the cost function given by
\begin{equation} \label{MScost}
C_{_{MS}}(a) = - \sum_{k=0}^{n-1} \log \Biggl( \int \exp \biggl( - \frac{d_k^2(x)}{2\epsilon^2} \biggr) \mu_a(dx) \Biggr),
\end{equation}
where $d_k^2(x) =  |(y_k,y_{k+1}) - (x,T_a(x))|^2$, $\epsilon$ is the variance of the noise process $(\epsilon_n)_n$, and $\mu_a$ is a particular invariant measure for the map $T_a$. They argue that the minimum of $C_{_{MS}}$ provides more reliable parameter estimates due to its inclusion of information regarding the invariant measure $\mu_a$. It is perhaps a shortcoming of this method that one must know the variance of the noise process and the invariant measure $\mu_a$ in order to calculate $C_{_{MS}}(a)$. In practice, the authors suggest approximating the integral with respect to $\mu_a$ by a sum over a long piece of trajectory simulated from the model in the hopes that this approximation will be close to the integral by the ergodic theorem. The authors provide numerical evidence that $C_{_{MS}}$ provides better parameter estimates than either $C_{_{OSLS}}$ or $C_{_{TLS}}$, although again no theoretical results are available to justify this comparison.

Meyer and Christensen \cite{MeyerChristensen2000}, following up on the work of McSharry and Smith \cite{McSharrySmith1999}, propose to model the system using a combined noise state-space model of the form (\ref{generalobs})-(\ref{generaldyn}), and proceed via an MCMC algorithm for performing the inference. In particular, they take a Bayesian approach, modeling both the true states $X_n$ and the parameters $a$ as unknown variables. They assume that the process $(X_n)_n$ forms a Markov chain (by adding dynamical noise to the model). Then they compute posterior probabilities of the unobserved variables using the Gibbs sampler and the Metropolis-Hastings algorithm.

In his paper \cite{Judd2003_2}, titled ``Chaotic-time-series reconstruction by the Bayesian paradigm: Right results by wrong methods," Judd discusses the Bayesian approach of Meyer and Christensen \cite{MeyerChristensen2000}, and argues that their approach might work, but for ``accidental" reasons. In particular, he objects to the fact that Meyer and Christensen have replaced the deterministic model by a stochastic model, and he claims to formulate the ``correct" Bayesian approach for the deterministic model, which he acknowledges is essentially that given by Berliner \cite{Berliner1991} (presented in Section \ref{SectionIntroToLikelihoods}). He argues that their model only appears to give correct results because it happens to find shadowing trajectories of the true system. (We remark that an $\epsilon$ shadowing trajectory for a sequence $(z_n)_n$ of states in $\cX$ is a true orbit $(x_n)_n$ of the system such that $d(x_n,z_n) < \epsilon$ for all $n$.)  In the end, he argues for methods based on a direct search for shadowing trajectories, which purportedly require significantly less computational effort than the Bayesian approach of Meyer and Christensen. Such methods are often referred to as gradient descent methods; for further reading about these methods, see \cite{Judd2008,RidoutJudd2002} and references therein.

Pisarenko and Sornette \cite{PS2004} consider the parameter estimation methods discussed above, as well as the method of moments. They point out that the method of moments seems to be the only method so far considered whose asymptotic consistency has been rigorously proved on even a single example. They also provide a careful analysis of the work of McSharry and Smith. This analysis sheds light on some errors, both quantitative and qualitative, in the work \cite{McSharrySmith1999}. In order to provide a useful estimation procedure, they propose a ``pure" likelihood method, in which they cut the time-series data into $n_1$ sub-intervals of length $n_2$ and perform ML estimation on each interval independently. In this method, the resulting $n_1$ parameter estimates are averaged to produce a single estimate. The motivation behind their method seems to be the following. A theoretically true/pure ML method involves treating $x_0$ as a parameter to be estimated (as in Section \ref{SectionIntroToLikelihoods}), but the chaotic nature of the system means that the system forgets its initial condition exponentially quickly, which implies that it cannot be reliably estimated. Hence, they arrive at the method of chopping the time series into smaller pieces (which hopefully still contain useful information about the initial condition of each piece) and using the pure ML method on each piece. The statistician interested in mathematical rigor is likely to find this work rewarding to read. A word of caution: they conclude their article by stating that ``the situation is rather hopeless for the establishment of a meaningful statistical theory of estimation using the continuous theory of classical statistics to such discontinuous objects as the invariant measures of chaotic dynamical systems."

Smirnov \textit{et al} \cite{SVP2005} note that the piecewise ML method of Pisarenko and Sornette \cite{PS2004} suffers from significant bias and potentially large variance, since it relies on chopping the data into many small subsets. In the case of one-dimensional maps, the authors propose a method based on backwards iteration of the map. They interpret their method as also relying on a ML principle, and they claim (with numerical support but no proof) that their method is asymptotically consistent with variance typically decreasing like $n^{-2}$, where $n$ is the length of the observed time series.

The work of Horbelt and Timmer \cite{HT2003} seeks to quantify the rate of convergence of parameter estimates to the true parameter value in the observational noise case as the number of observations grows. In the introduction, the authors claim that the MLE in this setting is unbiased and efficient, for which they refer the reader to an earlier version of the book \textit{Theory of point estimation} by Lehmann and Casella \cite{Lehmann1998}, although it seems clear that this statement is mistaken, since in some cases it can be shown to be asymptotically biased.  Nonetheless, the authors find numerical evidence for various scaling laws of the variance of the MLE.

The work of Nakamura \textit{et al} proposes yet another parameter estimation method in the observational noise setting \cite{NHJKS2007}. Here the authors suggest an iterative method that alternates between estimating the system states and the system parameters. In each of the optimization steps, they use somewhat standard techniques. To estimate the parameters, they minimize the cost function $C_{_{OSLS}}$ in (\ref{OSLScost}) with respect to the parameter $a$. To estimate the states, they state that one could use any filtering method, such as the extended Kalman filter \cite{WalkerMees1997,WalkerMees1998} or gradient descent noise reduction (also known as gradient descent state estimation) \cite{JuddSmith2001,KostelichSchreiber1993,RidoutJudd2002}. The novelty of their approach lies in the fact that they iterate between these two estimation steps.

\subsection{Method of moments} \label{SubSectionMethodOfMoments}

Here we mention a method of parameter estimation that has been shown to be consistent at least for the logistic family, discussed in Example \ref{ExampleLogistic}. For the observational noise model, this method, discussed in \cite{PS2004}, appears to be the only method that has been proved to be consistent for at least one non-trivial example.

We consider the model (\ref{obsnoiseobs})-(\ref{obsnoisedyn}), where $\cX = [-1,1]$, $\cA = [0,2]$, and $T_a(x) = 1-ax^2$, which is change of coordinates of the family in Example \ref{ExampleLogistic}. Assume that the underlying trajectory process $(X_n)_n$ is ergodic, which is the case if one assumes that $X_0$ is drawn from an ergodic invariant measure $\mu_a$ for the map $T_a$. Alternatively, one may assume that $a$ is chosen such that $T_a$ has an acim $\mu_a$ (as discussed in Example \ref{ExampleLogistic}) and $X_0$ is drawn from Lebesgue measure. Also assume that the observational noise is additive (\textit{i.e.} $Y_n = X_n + \epsilon_n$) and $(\epsilon_n)_n$ is i.i.d. Gaussian with mean $0$ and variance $\epsilon^2$. For a sequence $(z_k)_{k=0}^{n-1}$, let $A_n(z_k) = \frac{1}{n} \sum_{k=0}^{n-1} z_k$ and for any $f : \R \to \R$, let $\EE_{\mu_a}(f) = \int f(x) \, d\mu_a(x)$.  Then by the ergodic theorem
\begin{align}
\lim_{n \to \infty} A_n(Y_k) & = \EE_{\mu_a}(x)  \label{EqnMethodOfMoments1} \\
\lim_{n \to \infty} A_n(Y_k^2) & = \EE_{\mu_a}(x^2)  \label{EqnMethodOfMoments2} \\ 
\lim_{n \to \infty} A_n(Y_k^3) & = \EE_{\mu_a}(x^3) + 3 \epsilon^2 \EE_{\mu_a}(x) \label{EqnMethodOfMoments3} \\  
\lim_{n \to \infty} A_n(Y_kY_{k+1}) & = \EE_{\mu_a}(x) - a \EE_{\mu_a}(x^3).  \label{EqnMethodOfMoments4}
\end{align}
Also, averaging the equation $x_{n+1} = 1 - a x_n^2$, we obtain that
\begin{equation}
\EE_{\mu_a}(x) = 1 - a \EE_{\mu_a}(x^2).  \label{EqnMethodOfMoments5}
\end{equation}
Combining Equations (\ref{EqnMethodOfMoments1})-(\ref{EqnMethodOfMoments5}), we arrive at the following estimates for the unknown parameters $a$, $\EE_{\mu_a}(x)$, $\EE_{\mu_a}(x^2)$, $\EE_{\mu_a}(x^3)$, and $\epsilon$:
 \begin{align*}
 \hat{a}_n & = \frac{A_n(Y_k Y_{k+1}) + 2 A_n(Y_k) + 3 \bigl(A_n(Y_k)\bigr)^2 }{ 3 A_n(Y_k) \bigl(A_n(Y_k)\bigr)^2 - A_n(Y_k^3)} \\
 \EE_{\mu_a}(\hat{x})_n & = A_n(Y_k) \\
 \EE_{\mu_a}(\hat{x}^2)_n & = A_n(Y_k^2) - \hat{\epsilon}_n \\
 \EE_{\mu_a}(\hat{x}^3)_n & = \frac{1}{\hat{a}_n} \bigl( A_n(Y_k) - A_n(Y_k Y_{k+1})\bigr) \\
 \hat{\epsilon}_n & = \frac{ A_n(Y_k^3) - \EE_{\mu_a}(\hat{x}^3)_n }{3 A_n(Y_k)}
 \end{align*}
These estimates are consistent by the ergodic theorem, but they might converge quite slowly, as there is no general rate of convergence in the ergodic theorem.
%%%%%%%%%%%%%%%%%%%%%%%%%%%%%%%%%%%%%%%%%%%%%%%%%%%%%%%%%%%%%%%%%%%%
%%%%%%%%%%%%%%%%%%%%%%%%%%%%%%%%%%%%%%%%%%%%%%%%%%%%%%%%%%%%%%%%%%%%
\section{Dynamical noise only} \label{SectionDynNoise}

If only dynamical noise is present in the model (\ref{generalobs})-(\ref{generaldyn}), then the system (\ref{generalobs})-(\ref{generaldyn}) reduces to the following situation:
\begin{align}
Y_n & = f_a(X_n) \label{dynnoiseobs} \\
X_{n+1} & = T_a(X_n,\delta_n), \label{dynnoisedyn}
\end{align}
where $(\delta_n)_n$ is a noise process, $T : \cA \times \cX \times \cN \to \cX$ is a parametrized family of noisy maps, and $f : \cA \times \cX \times \cN \to \cY$ is a parametrized family of observation functions. The dynamical noise model has been studied in the dynamical systems literature under the name ``random dynamical systems'' (see \cite{KiferLiu2006} and references therein). The process $(X_n)_n$ forms a discrete-time Markov chain on the continuous state space $\cX$ (see the book of Meyn and Tweedie \cite{MeynTweedie2009} and references therein). In this case, some of the estimation methods from the statistical literature on time series and state space models may apply.

Without using this Markov structure, Adams and Nobel have studied the non-parametric reconstruction of such systems from direct observations (\textit{i.e.} $Y_n = X_n$) \cite{Nobel2001,AdamsNobel2001_2}. In particular, they used adaptive histogram methods to show results similar to those regarding non-parametric reconstructions of systems with no noise, as in Section \ref{SectionNoNoiseReconstructionDirectObservations}. These methods do not work in the observational noise case precisely because in that setting they suffer from the problem of errors in variables, as discussed in Section \ref{SubSectionVariations}.

A common setting for random dynamical systems is to assume that there is a map $T : \cX \to \cX$, where $\cX$ is a compact manifold and $T$ is smooth, with a ``natural'' invariant probability measure $\mu$. In common examples, $T$ might be a (non-uniformly) hyperbolic map and $\mu$ might have the property that almost every initial condition with respect a volume measure on the manifold equidistributes with respect to $\mu$. In such cases, one typically adds dynamical noise as follows. Let $\epsilon >0$. For each $x$ in $\cX$, let $\mathbb{P}_{\epsilon}(x,\cdot)$ be the uniform measure on the ball of radius $\epsilon$ about the point $T(x)$. Then the Markov chain corresponding to this random dynamical system is determined by viewing $\mathbb{P}_{\epsilon}$ as the transition kernel for the chain. Under some conditions, the chain corresponding to $\mathbb{P}_{\epsilon}$ will have a unique stationary distribution, $\mu_{\epsilon}$. A well-known result (see \cite{KiferLiu2006}) states that under certain conditions, the measure $\mu_{\epsilon}$ converges to $\mu$ weakly as $\epsilon$ tends to $0$. To the best of our knowledge, no theoretical work on parameter estimation has been conducted for this particular setting, perhaps making it an area ripe for progress. On the other hand, this setting may be viewed as a particularly degenerate version of the general state-space setting, in which there is no observational noise, and therefore all methods described in Section \ref{SectionCombNoise} may also be applied here.

%%%%%%%%%%%%%%%%%%%%%%%%%%%%%%%%%%%%%%%%%%%%%%%%%%%%%%%%%%%%%%%%%%%%
%%%%%%%%%%%%%%%%%%%%%%%%%%%%%%%%%%%%%%%%%%%%%%%%%%%%%%%%%%%%%%%%%%%%
\section{General state space models} \label{SectionCombNoise}

In this section we consider the full system (\ref{generalobs})-(\ref{generaldyn}), where both dynamical noise and observational noise are present. Specific versions of such models have long been considered in the statistics literature, where they are known as state space models \cite{DurbinKoopman2001}. The literature on state space models in both applied and theoretical statistics is extensive and  \cite{Hamilton1994,Pole1994} are two excellent texts covering applied modeling on this topic. The models can be summarized as the study of hidden Markov models (HMMs) in general state-spaces. (For an article discussing the connections between ergodic theory and finite state HMMs, see \cite{BoylePetersen2011}.) Theoretical understanding of general HMMs has been a challenge and rigorous statements on consistency in parameter estimation have only appeared recently \cite{DMOvH2011} (see Section \ref{SectionMLEforHMMs}). Most of the work in this area has been devoted to the problem of state estimation or filtering, and even at a computational level the problem of parameter estimation is still largely unsolved. In this section we survey some of the most studied  approaches to filtering and discuss parameter estimation where there are results.

%By including the unknown parameters in a new state space, one may trivially transfer any filtering method into a parameter estimation. Unfortunately, by artifically viewing the parameters as part of the state space, one may lose much of the structure of the parameter inference problem, and furthermore one often introduces significant degeneracy into the model, which may severely damage the performance of the filtering method.
\subsection{Kalman filter and some generalizations}
The simplest such models assume that the dynamics are linear and the noise is additive Gaussian:
\begin{align*} \label{EqnLinearGaussian}
X_{n+1} & = A X_n + B \delta_{n+1} \\
Y_n & = C X_n + D \epsilon_n,
\end{align*}
where here $A$, $B$, $C$, and $D$ are all matrices of the appropriate dimension and $(\delta_n)_n$ and $(\epsilon)_n$ are independent i.i.d. Gaussian processes. In this case, the optimal solution to the state estimation or denoising problem is given by the well-known Kalman filter \cite{Eubank2006,Kalman1960}. Generalizations of the ideas behind Kalman filtering to non-parametric models have been an extensive area of research in Bayesian and frequentist inference \cite{DurbinKoopman2001,FSJW2010,FSJW2011,Muhkerjee2011}. 

Conceptually, the simplest generalization of the Kalman filter to nonlinear models involves linearizing the models at each time point and then using the Kalman filter. This method is often called the \textit{extended Kalman filter} (EKF) \cite{Jazwinski1970,AndersonMoore1979}. While the Kalman filter is optimal in the sense that is the minimal-variance unbiased estimator, the general EKF is known to be biased. Furthermore, due to the linearization of the model, the propagation of the error covariance estimates may behave quite poorly if the non-linear terms in the model are significant.

The unscented Kalman filter (UKF) \cite{JulierUhlmann1996,JulierUhlmann1997} provides a deterministic sampling scheme that has been observed to outperform the EKF. The basic idea behind the UKF is that instead of approximating the model by linearization, one ought to use the exact model but approximate the posterior distributions by Gaussian distributions. The sampling scheme is designed to insure that the first two moments of the posterior distributions match the first two moments of the approximating distributions. It is believed that the UKF outperforms the EKF because it may be viewed as an unbiased second-order method, whereas the EKF is a biased first-order method. Of course, the UKF is believed to have shortcomings of its own; in particular, it assumes that the posterior distributions are Gaussian, which is certainly not the case in general. Also, the number of samples required for the UKF is at least the dimension of the state space, and in high-dimensional settings this fact makes the UKF computationally intractable. A wide variety of Monte Carlo (MC) methods have been proposed to overcome these issues.

Another generalization of the Kalman filter is known as the ensemble Kalman filter (EnKF) \cite{Evensen1994,BLE1998,Evensen2003}. This method is a Monte Carlo method that is particularly popular in the weather prediction community. In fact, this method may be thought of as a type of particle filter (see Section \ref{sec:smc}).

\subsection{MLE for HMMs} \label{SectionMLEforHMMs}

If one is willing to consider point estimates of unknown parameters in a setting where the likelihood function is known, then one can consider the maximum likelihood method (MLE) for parameter estimation. Let us now state the main result of the paper \cite{DMOvH2011}, which gives sufficient conditions for the consistency of MLE in this context. Let $(X_k,Y_k)_{k=1}^{\infty}$ be a hidden Markov model (HMM) of the form (\ref{generalobs})-(\ref{generaldyn}). Let $a^*$ denote a fixed parameter value  in $\cA$. Assume that the HMM with parameter $a^*$ has a unique stationary distribution, and let $\mathbb{P}_{a^*}$ be the corresponding stationary HMM. Denote by $p^{\nu}(y_0^n,a)$ the likelihood of the observations $Y_0^n$ with initial distribution $X_0 \sim \nu$ and parameter $a$. Consistency of the maximum likelihood estimator (MLE) may now be stated in the following form: if $a_n = \argmax_a p^{\nu}(y_0^n,a)$, then $a_n$ converges $\mathbb{P}_{a^*}$-a.s. to $a^*$ as $n$ tends to infinity. The main result of \cite{DMOvH2011} gives some general conditions under which the MLE is consistent in this sense. A precise statement of these general conditions is beyond the scope of this survey.

\subsection{Bayesian inference} \label{SectionBayesianInference}

Recall the Bayesian formulation of state space estimation or filtering. Here one assumes that the model (\ref{generalobs})-(\ref{generaldyn}) gives rise to probability densities $\mu(x_0)$, $p(x|x')$, and $q(y|x)$, which define the initial distribution, transition kernel, and marginal distribution of the observation process, respectively. The densities are with respect to some fixed reference measures denoted $dx$ and $dy$. In this framework, we are given access to finitely many observations $y_0^{n-1}$, and we would like to estimate the true trajectory $x_0^{n-1}$. Our assumptions define likelihood functions
\begin{equation*}
p(x_0^{n-1}) = \mu(x_0) \prod_{k=0}^{n-2} p(x_{k+1}|x_k),
\end{equation*}
and 
\begin{equation*}
p(y_0^{n-1}|x_0^{n-1}) = \prod_{k=0}^{n-1} q(y_k|x_k).
\end{equation*}
Given the observations $y_0^{n-1}$, the posterior distribution for $X_0^{n-1}$ is given by
\begin{equation*}
p(x_0^{n-1}|y_0^{n-1}) = \frac{p(x_0^{n-1},y_0^{n-1})}{p(y_0^{n-1})},
\end{equation*}
where
\begin{align*}
& p(x_0^{n-1},y_0^{n-1}) = p(x_0^{n-1}) p(y_0^{n-1} | x_0^{n-1})\\
& p(y_0^{n-1}) = \int p(x_0^{n-1},y_0^{n-1}) dx_0^{n-1}.
\end{align*}

There are a few instances when these distributions may be calculated analytically, such as when the system is linear and the noise is Gaussian or when $\{X_n\}_n$ is a finite state Markov chain. Outside of these cases, there is no analytical method for calculating the posterior distribution, and therefore one seeks a numerical approximation for this distribution. With the significant advances in computational power in recent years, there has been a remarkable amount of research devoted to finding efficient computational approaches to approximating such posterior distributions. In the remainder of this section we briefly discuss the linear Gaussian case and a few of its generalizations to nonlinear or non-Gaussian situations. Section \ref{SectionMC} discusses some of the more recent computational approaches to filtering.

An interesting work in the Bayesian context is \cite{Shal2009} where the author studies posterior consistency for dependent data from an information theoretic point of view.  The author establishes posterior consistency  for misspecified models under the assumption of asymptotic equipartition property. For finite state space ergodic models, this is implied by the Shannon-McMillan-Breiman theorem. It could be interesting and useful to extend the ideas from \cite{Shal2009} to prove posterior consistency in parameter estimation for more general dynamical systems.

\subsection{Inference for dynamical systems via simulation based methods} \label{SectionMC}

In the general non-linear, non-Gaussian state-space setting of (\ref{generalobs})-(\ref{generaldyn}), the posterior distributions for $x_0^{n-1}$ are not available in closed form, as they involve some integrals for which no analytical evaluation methods exist. In order to perform inference in this setting, a great deal of effort has been devoted to developing sophisticated computational algorithms for for sampling from these posterior distributions. One general idea is to use Monte Carlo (MC) methods to estimate the integrals of interest. It is worth emphasizing that there has been a huge amount of work in this direction, and we do not claim to provide a comprehensive survey of all the relevant results. For an introduction to MC methods, see the book \cite{RobertCasella2004}.

\subsubsection{MCMC methods, SMC and Particle Filters.} \label{sec:smc}

If one cannot sample from the posterior distribution directly, then one often turns to Markov chain Monte Carlo (MCMC) methods. For a discussion of such methods, see the books \cite{RobertCasella2004,WestHarrison1997} and references therein. 
%Here the idea is to construct an ergodic Markov chain whose invariant distribution is the desired posterior distribution. Then, by the ergodic theorem, a typical path of the chain may be used to estimate basic properties of the posterior distribution, such as mean or variance. 
Such methods have been used for  parameter estimation in dynamical systems 
(\textit{e.g.}, \cite{CMKL2001}).

Traditional Monte Carlo or MCMC methods may be used to perform ``batch" inference, \textit{i.e.} when all of the observations are available at once and one would like to estimate $p(x_0^{n-1}|y_0^{n-1})$ for fixed $n$, although even in this setting they might be prohibitively computationally expensive. When the goal is to perform ``on-line" or sequential inference, or in an effort to try to reduce the computational expense, one might try sequential Monte Carlo methods (SMC) and their many variations. A particularly popular version of these methods is known as particle filtering. For a well-written, thorough introduction to the principles of sequential Monte Carlo (SMC) and particle filtering methods, see the recent tutorial by Doucet and Johansen \cite{DoucetJohansen2011}. For an incomplete list of works concerning SMC and particle filtering, as well as their adaptations to parameter estimation, see \cite{BLE1998,DelMoral2004,dMDS2010,Evensen1994,Evensen2003,Fearnhead2002,IBAK2011,KDSM2009,LopesTsay2011,MukherjeeWest2009,OCDM2008,PDS2011,Storvik2002}. The basic idea is that the posterior distributions of interest are approximated by a finite collection of $N$ samples, called particles, which are recursively propagated through the model. The main theoretical advantage of these methods is that one is often able to establish the convergence of the approximations to the true posterior distributions as the number of particles $N$ tends to infinity.

% Is EnKF a particle filtering method? One MC method popular in the weather prediction community is known as the Ensemble Kalman filter (EnKF) \cite{Evensen1994,BLE1998,Evensen2003}.

%\cite{MukherjeeWest2009}

%Overview of param. estimation for HMM: \cite{KDSM2009} 

%Mention degeneracy issue and attempts to mitigate?

%Reducing the Variance using Rao-Blackwellized Particle Filters?

%SMC Smoothing?

%SMC for on-line Bayesian static parameter estimation in state-space models
%\cite{Fearnhead2002,Storvik2002,LopesTsay2011}

%A pragmatic approach consists of adding an artificial dynamic noise on the static parameter (references to include.)

%SMC for on-line and batch maximum likelihood inference of static parameter estimation in state-space models 
%\cite{OCDM2008,dMDS2010,IBAK2011} \cite{PDS2011}

%SMC for batch Bayesian static parameter estimation in state-space models
%Obviously for batch joint Bayesian state and parameter estimation, you can use MCMC methods. However it can be difficult to design efficient algorithms. In the context where one can only simulate the latent process but does not have access to the transition prior, standard MCMC just fail. SMC can come to the rescue in these scenarios.

%SMC as an alternative/complement to MCMC

\subsubsection{ABC methods.}

Most of the methods mentioned previously in this section rely on explicit knowledge and evaluation of the likelihood function. In many situations, such as in high dimensional complex models, the likelihood function may not be available or is 
computationally expensive to evaluate.  In such scenarios, a simple computational method called approximate Bayesian computation (ABC)  offers a powerful alternative to conduct statistical inference. ABC was first proposed as a philosophical argument in \cite{Rubi1984} and introduced to population genetics  in \cite{TBGD1997}. Since then these methods have become extremely popular in many applied fields. A partial list of references include \cite{DJSP2011,dMDJ2011,McKinleyCookDeardon2009,PWS2010,PSPF1999,RAWS2009,SFT2007,TWSIS2009,Wilkinson2008}. A  good review with applications to filtering 
is \cite{Mari2011}. Briefly speaking,  in ABC methods one first draws a parameter value $\theta^*$ from the prior distribution and generates synthetic data from the likelihood model corresponding to $\theta^*$. If the synthetic data   ``is similar to" the observed data (measured in some metric) up to a prespecified tolerance then $\theta^*$ is accepted as a draw from the (approximate) posterior distribution. Choosing the metric and the tolerance level are difficult problems, but partial results are known (\cite{FearPrang2012}).

  An  important point to note is that in many examples,  a summary statistic instead of the original data set  is used for matching. This clearly results in loss of information (and sometimes even results in invalid inference; see \cite{Robetal2011}) and thus raises the interesting question about when one can perform consistent model selection using the ABC methodology. In \cite{Mari2011:02} a sufficient criteria is worked out, but clearly more needs to be done especially in the context of dynamical systems.

%%%%%%%%%%%%%%%%%%%%%%%%%%%%%%%%%%%%%%%%%%%%%%%%%%%%%%%%%%%%%%%%%%%%
%%%%%%%%%%%%%%%%%%%%%%%%%%%%%%%%%%%%%%%%%%%%%%%%%%%%%%%%%%%%%%%%%%%%
\section{Open questions and future directions} \label{SectionOpenQuestions}

Here we list some open questions related to parameter inference in dynamical systems and discuss possible future research directions.

The first question examines if parameter estimation is possible even if denoising is impossible. 
\begin{ques} \label{QuestionLalleyNobel}
 As shown by Lalley and Nobel \cite{LalleyNobel2006}, there are instances in which state estimation or denoising in the observational noise setting is impossible. Is it possible to exhibit a family of topological dynamical systems $(X,T_a)$ on a compact metric space $X$ such that consistent denoising is (provably) impossible but consistent parameter estimation is nonetheless (provably) possible?
\end{ques}

The most common method in theory and practice for parameter estimation is ML. It is open if ML in the observational noise setting is consistent. A related question
is can the approach to proving consistency results for HMMs in  \cite{DMOvH2011} be adapted to the observational noise setting.
\begin{ques} \label{QuestionMLconsistency}
 Recall the definition of the ML estimator of the parameter $a$ given in (\ref{EquationMLE}):
\begin{equation*}
\hat{a}_n = \underset{a}{\argmax} \max_{x_0} L(x_0,a),
\end{equation*}
where $L(x_0,a)$ is the likelihood of $x_0$ and $a$ conditional on the observations $(y_k)_{k=0}^{n-1}$. In the observational noise case setting (\ref{obsnoiseobs})-(\ref{obsnoisedyn}), what are necessary and sufficient conditions on the system such that $\hat{a}_n$ converges to $a$ with probability $1$ (for almost every initial condition $x_0$ with respect to an ergodic measure $\mu$)? If necessary and sufficient conditions are out of reach given current tools, partial answers to this question in the form of general sufficient conditions might also be interesting.
\end{ques}

In order to get finite sample error bounds, one would also like to know about the deviations of the MLE from its average. This line of reasoning leads to the following question.
\begin{ques} \label{QuestionMLdeviations}
In the observational noise setting (\ref{obsnoiseobs})-(\ref{obsnoisedyn}), under which conditions on the system is it true that the MLE is asymptotically normal in the observational noise setting?
\end{ques}

Since the method of moments presented in Section \ref{SubSectionMethodOfMoments} is currently the only example to our knowledge for which consistency of any parameter estimation method can be proved in the observational noise setting, it is worth considering how it might be generalized.
\begin{ques} \label{QuestionMethodOfMoments}
 Can the method of moments presented in Section \ref{SubSectionMethodOfMoments} for the logistic family be generalized? Under what conditions is it applicable and consistent?
\end{ques}

In the combined noise setting of Section \ref{SectionCombNoise}, it is still the case that the issue of parameter inference has not been satisfactorily resolved. Certainly any filtering method may be trivially extended to a parameter estimation algorithm by extending the state space to include the parameters, but in such cases the degeneracy of the extended system typically causing the filtering methods to fail. Let us paraphrase a question in \cite{DoucetJohansen2011}.

\begin{ques} \label{QuestionCombNoise}
 Under what conditions on the model are there efficient algorithms for parameter estimation in the general state space setting? What theoretical guarantees can be given to justify such algorithms?
\end{ques}

The range of applications of statistical inference methods for deterministic dynamical systems seems to be increasing rapidly. These systems present significant new challenges, since the deterministic systems may have very long-range dependency structures. It would be a significant breakthrough if methods could be developed that provided asymptotically consistent algorithms for parameter estimation; moreover, one would like to have finite-size sample bounds on the accuracy of these algorithms. Given the difficulty of dealing with the long-range dependencies present in general in the observational noise model, it appears likely that the traditional methods of parameter inference may not work particularly well in this setting, and therefore new ideas and methods should be developed.

One possible approach would be to consider a weakened notion of consistency.  For example, one could consider a parameter estimation method to be consistent if it returns a set of plausible parameters that asymptotically contains the true parameter. Such weakened notions of consistency might be necessary for providing some theoretical justification of parameter estimation algorithms when achieving strong consistency appears out of reach.

Let us close with one recent development in the field of dynamical systems and ergodic theory that might be useful in obtaining such rates of convergence. The concentration inequalities mentioned at the end of Section \ref{SectionErgTheory} provide a powerful method for obtaining finite sample error bounds for a wide class of statistical estimators for a wide class of dynamical systems. One might hope that these concentration inequalities can be used to get rigorous error bounds for parameter estimation algorithms.

\vspace{.2in}

\noindent \paragraph{\bf Acknowledgments} The authors would like to thank Andrew Nobel, Ramon van Handel, John Harer,
Konstantin Mischaikow, Christian Robert, Mark Girolami and Andrew Stuart for discussions, comments and help with references. SM and KM would like to acknowledge AFOSR FA9550-10-1-0436 and NSF DMS-1045153
for partial support. SM would also like to acknowledge NSF CCF-1049290 for partial support. NSP would like to thank NSF for partial support through the grant NSF DMS-1107070.

%%%%%%%%%%%%%%%%%%%%%%%%%%%%%%%%%%%%%%%%%%%%%%%%%%%%%%%%%%%%%%%%%%%%%

\bibliographystyle{plain}
\bibliography{SurveyRefs}
%\bibliography{InferenceRefs,ConleyIndexRefs}

%\nocite{*}

\end{document}